\documentclass[a4paper,12pt]{amsart}

\usepackage{fullpage}

\makeindex

\usepackage{xcolor}
\usepackage[
colorlinks=true,
citecolor=blue,
urlcolor=blue,
linkcolor=blue,
pdfborder={0 0 0},
breaklinks
]{hyperref}

\usepackage{amsfonts,graphics,amsmath,amsthm,amscd,amssymb,amsmath,latexsym,euscript, enumerate,kotex,mathtools}
\usepackage{epsfig,url}
\usepackage{flafter}
\usepackage[all,cmtip,line]{xy}
\usepackage{array}
\usepackage[english]{babel}
\usepackage{overpic}
\usepackage{subfig}
\usepackage{multirow}
\usepackage{tikz-cd}
\usepackage[T1]{fontenc}
\usepackage{lmodern}
\usepackage{pifont}

\usepackage{microtype}
\usepackage{wrapfig}
\usepackage[shortlabels]{enumitem}
\setlist[enumerate, 1]{1\textsuperscript{o}}

\newtheorem{theorem}{Theorem}[section]
\newtheorem{lemma}[theorem]{Lemma}
\newtheorem{proposition}[theorem]{Proposition}

\newtheorem{corollary}[theorem]{Corollary}

\theoremstyle{definition} 
\newtheorem{definition}[theorem]{Definition}
\newtheorem{definition-lemma}[theorem]{Definition-Lemma}

\theoremstyle{remark}
\newtheorem{remark}[theorem]{Remark}

\numberwithin{equation}{section}

\newcommand{\Q}{\mathbb{Q}}

\DeclareRobustCommand{\O}{\mathcal{O}}

\newcommand{\Ker}{\mathrm{Ker}}

\DeclareMathOperator{\Hom}{\mathrm{Hom}}
\DeclareMathOperator{\Ext}{\mathrm{Ext}}

\def\Spec{\operatorname{Spec}}

\def\codim{\operatorname{codim}}

\let\oldframe\frame
\renewcommand\frame[1][allowframebreaks]{\oldframe[#1]}

\title[On a cohomological property of the center of a resolution]
{On a cohomological property of the center of a resolution}

\begin{document}

\author[D.~Kim]{Donghyeon Kim}
\address[Donghyeon Kim]{Department of Mathematics, Yonsei University, 50 Yonsei-ro, Seodaemun-gu, Seoul 03722, Republic of Korea}
\email{narimial0@gmail.com}

\date{\today}

\begin{abstract}
In this note, we explore the cohomological property of the codimension of the center of a resolution. In particular, we define a resolution $f:X'\to X$ to be $q$-birational if the center of $f$ satisfies $\codim \mathrm{Cent}(f)\ge q+1$, and we prove that $R^if_*\mathcal{O}_X(E)=0$ for every $1\le i\le q-1$ and every $f$-anti-nef effective $f$-exceptional divisor $E$ on $X'$ if $X$ is $(R_q)$ and $(S_{q+1})$. We also discuss a partial converse of the theorem.
\end{abstract}

\maketitle

\section{Introduction}
Throughout the note, $k$ is an algebraically closed field of characteristic $0$, and all varieties are separated, finite type, reduced, and irreducible schemes over $k$. Let us assume $q \ge 1$ is a positive integer.

\smallskip

Let $X$ and $X'$ be normal projective varieties, and let $f \colon X' \to X$ be a proper birational morphism. Assume there is an effective $f$-exceptional divisor $F$ on $X'$. Then the following statement is standard and almost immediate:
$$ f_*\mathcal{O}_{X'}(F) = \mathcal{O}_X. $$
A central goal of this note is to generalize this equality in the following way.

\begin{definition}
Let $X,X'$ be any smooth varieties and $f:X'\to X$ be any proper birational morphism. We say that $f$ is \emph{$q$-birational} if the codimension of the center of $f$ is $\ge q+1$.
\end{definition}

\begin{theorem}\label{-1application}
Let $X,X'$ be any projective varieties, $X$ be $(R_q)$ and $(S_{q+1})$, let $X'$ be smooth, and let $f:X'\to X$ be any proper birational morphism. If $f$ is a $q$-birational morphism, then $R^if_*\O_{X'}(E)=0$ for every $1\le i\le q-1$, and every $f$-anti-nef effective $f$-exceptional divisor $E$ on $X'$.

\smallskip

Conversely, if $X$ is Cohen--Macaulay, and if there is an $f$-anti-ample effective $f$-exceptional divisor $E$ on $X'$ such that $R^if_*\O_{X'}(tE)=0$ for every $1\le i\le q-1$ and every $t\gg 0$, then $f$ is $q$-birational.
\end{theorem}

\smallskip

Note that if $X$ is $\Q$-factorial, then there always exists an $f$-anti-ample effective $f$-exceptional divisor on $X$ (cf. \cite[Lemma 3.6.2 (3)]{BCHM10})

\smallskip

The idea of the proof is the use of the Kawamata–Viehweg vanishing theorem with local cohomology computation. The proof of the partial converse is inspired by \cite[The proof of Lemma 3.3]{Kovac}.

\section{Preliminaries}
The aim of the section is to collect properties and theorems which will be used in the next sections. From this section and until the end of the note, we assume $2\le q<\dim X$ and for any variety $X$ with a point $x\in X$, $\dim x$ is the dimension of the closure of $x$, $\codim_X x:\dim X-\dim x$, and $X_x:=\Spec \O_{X,x}$.

\smallskip

For a proper morphism $f:X'\to X$ between two varieties $X'_x:=X'\times_X X_x$, let us denote the structural morphism by $f:X'_x\to X_x$. Moreover, $\O_{X',x}:=f^*\O_{X,x}$ is the structure sheaf of $X'_x$, and $\omega_{X,x}$ is the dualizing sheaf of $X'_x$.

\smallskip

The following theorem is well-known \emph{relative Kawamata--Viehweg vanishing theorem}.

\begin{theorem}[{cf. \cite[Theorem 1-2-3]{Kawamata}}]
Let $X'$ be any smooth projective variety and $X$ be any projective variety. Let $f:X'\to X$ be any proper birational morphism. For any $f$-nef Cartier divisor $D$ on $X'$, $R^if_*\O_{X'}(K_{X'}+D)=0$ for $i\ge 1$.
\end{theorem}

Let us introduce the following concepts.

\begin{definition}
Let $X$ be any variety.
\begin{itemize}
    \item[(a)] We say $X$ is $(R_q)$ if the singular locus of $X$ has codimension $\ge q+1$.
    \item[(b)] We say $X$ is $(S_{q+1})$ if for any point $x\in X$, the local cohomology $H^i_x(X_x,\O_{X,x})$ is zero for each $i<\min\{q+1,\codim_X x\}$.
\end{itemize}
\end{definition}

Let us define the \emph{dualizing complex}.

\begin{definition}
Let $X$ be a noetherian scheme. We denote by $\omega^{\bullet}_X$ the $(\dim X)$-shift of the normalized dualizing complex of $X$ (for the definition of \emph{normalized dualizing complex}, see \cite[0AU5]{Stacks}).
\end{definition}

In this note, we will use derived category machinery. Hence, it is worth stating the \emph{Grothendieck duality}.

\begin{theorem}[{See \cite[0AU3 (4)]{Stacks}}]
Let $X,X'$ be any noetherian schemes and $f:X'\to X$ be any proper birational morphism. For any $K\in D(X')$,
$$ R\mathcal{H}om_{\O_X}(Rf_*K,\omega^{\bullet}_{X}) \cong Rf_*R\mathcal{H}om_{\O_{X'}}(K,\omega^{\bullet}_{X'})$$
in $D(X)$, where $D(X)$ means the bounded derived category of coherent sheaves on $X$.
\end{theorem}

\begin{corollary}[{cf. \cite[Lemma 4.22]{Fujino}}]\label{ssg}
Let $X,X'$ be any varieties and $X'$ be smooth. For any proper birational morphism $f:X'\to X$ and point $x\in X$ with $\codim_X x\ge c$, $H^i_{f^{-1}(x)}(X',\O_{X'})=0$ for $i<c$.
\end{corollary}

\begin{proof}
We may use Grothendieck duality. Indeed, if $E$ is the injective hull of the residue field of $\O_{X,x}$ as an $\O_{X,x}$-module, we have
$$ 
\begin{aligned}
R\Gamma_{f^{-1}(x)}(X'_x,\O_{X'_x})&=R\Gamma_{x}(X_x,Rf_*\O_{X'_x})
\\ &=\Hom_{\O_{X,x}}(R\mathcal{H}om_{\O_x}(Rf_*\O_{X'_x},\omega^{\bullet}_{X,x}),E)
\end{aligned}
$$
where we used the Leray spectral sequence for the first equality and local duality for the second equality (for the statement of local duality, see \cite[Lemma 0AAK]{Stacks}). Moreover,
$$
\begin{aligned}
R\mathcal{H}om_{\O_x}(Rf_*\O_{X'_x},\omega^{\bullet}_{X,x})&=Rf_*R\mathcal{H}om_{\O_{X',x}}(\O_{X'_x},\omega_{X',x})
\\ &=Rf_*\omega_{X',x}=\omega_{X',x}
\end{aligned}
$$
where we used Grothendieck duality for the first equality and the Kawamata-Viehweg vanishing theorem for the second equality. Hence, we have the assertion.
\end{proof}

\section{Proof of the main theorem}
In this section, to simplify notions, for any noetherian scheme $X$, any coherent sheaf $\mathcal{F}$ on $X$, and any set $F$ of $X$, we denote $H^i_F(X,\mathcal{F}):=H^i_{\overline{F}}(X,\mathcal{F})$.

\smallskip

Let us define the following notions.

\begin{definition}
Let $X,X'$ be any normal varieties and $f:X'\to X$ be any proper birational morphism.
\begin{itemize}
\item[(a)] The \emph{center} of $f$ is the reduced closed subscheme $C$ of $X$, which is the image of the exceptional locus along $f$.
\item[(b)] We say $f$ is a \emph{$q$-birational morphism} if the center of $f$ has codimension $\ge q+1$.
\end{itemize}
\end{definition}

Let us prove the following birational property of being $(S_{q+1})$. The proposition is not used in the paper.

\begin{proposition}\label{good}
Let $X$ be any normal variety which is $(R_q)$. Then the following are equivalent:
\begin{itemize}
    \item[\emph{(a)}] $X$ is $(S_{q+1})$.
    \item[\emph{(b)}] For any $q$-birational morphism $f:X'\to X$ with $X'$ smooth, $R^if_*\O_{X'}=0$ for $1\le i<q$.
\end{itemize}
\end{proposition}

\begin{proof}
For $\text{(a)}\implies \text{(b)}$, choose any resolution $f:X'\to X$ of $X$ which is an isomorphism outside of the singular locus of $X$. For any point $x\in X$ of $X$, there is a spectral sequence
\begin{equation}\label{1}
E^{st}_2=H^s_x(X_x,(R^tf_*\O_{X'})_x)\implies H^{s+t}_{f^{-1}(x)}(X'_x,\O_{X',x}).
\end{equation}
Suppose that $x\in C$ is a point in the center $C$ of $f$. Note that $\codim_X x\ge q+1$.

\smallskip

Suppose that $2\le n\le q$ and $R^if_*\O_{X'}=0$ for $1\le i<n-1$. Then, by spectral sequence argument,
$$ E^{0(n-1)}_2=\cdots=E^{0(n-1)}_{n}$$
and $E^{0(n-1)}_{n+1}\subseteq E^{n-1}_{\infty}=H^{n-1}_{f^{-1}(x)}(X'_x,\O_{X',x})=0$ if $n-1\le q$ (this is due to Corollary \ref{ssg}). Hence
$$ E^{0(n-1)}_{n+1}=\Ker (E^{0(n-1)}_{n}\to E^{n0}_{n})=0$$
and $H^0_x(X_x,(R^{n-1}f_*\O_{X'})_x)=E^{0(n-1)}_{n}\subseteq E^{n0}_{n}.$
Moreover, $E^{n0}_{n}=E^{n0}_2=H^n_x(X_x,\O_{X,x})=0$ for $n\le q$ (this is due to the $(S_{q+1})$ condition on $X$). 

\smallskip

This means that $H^0_x(X_x,(R^{n-1}f_*\O_{X'})_x)=0$ and $R^{n-1}f_*\O_{X'}$ does not have associated points as $x$. Certainly, any associated point of $R^{n-1}f_*\O_{X'}$ is in $C$. Thus, $R^{n-1}f_*\O_{X'}=0$.

\smallskip

For $\text{(b)}\implies \text{(a)}$, let us assume that for any $q$-birational morphism $f:X'\to X$ with $X'$ smooth, $R^if_*\O_{X'}=0$ for $1\le i<q$ holds. It suffices to prove that $X$ is $(S_{q+1})$.

\smallskip

Consider \eqref{1}. What we have to prove is that for any point $x$ of $X$, $H^i_x(X_x,\O_{X,x})=0$ for $i<\min\{q+1,\codim_X x\}$. Indeed, by the condition, $H^0_x(X_x,(R^tf_*\O_{X'})_x)=0$ for $0\le t<q$. Thus, by the spectral sequence argument, we obtain
$$ E^{i0}_2=E^{i0}_3=\cdots=E^{i0}_{\infty}\subseteq H^i_{f^{-1}(x)}(X'_x,\O_{X',x})$$
for any $i\le q$ and using Corollary \ref{ssg} gives us the assertion.

\smallskip

Let us note that the above proof is similar to \cite[the proof of Lemma 5.12]{KM}. Furthermore, for $q=2$, \cite[Proposition 4.27]{Fujino} becomes our proposition if we assume $X$ is $(R_2)$.
\end{proof}

\smallskip

Let us start proving the main result of this section, Theorem \ref{-1application}.

\begin{lemma}\label{zerothapplication}
If $X$ is a variety with $(R_q)$ and $(S_{q+1})$, then for any $q$-birational morphism $f:X'\to X$ from a smooth variety, $R^if_*\O_{X'}(E)=0$ for $1\le i<q$ and an anti-nef effective $f$-exceptional divisor $E$.
\end{lemma}

\begin{proof}
Let us note that $R^if_*\O_{X'}=0$ for $1\le i<q$ by Proposition \ref{good}.

\smallskip

We may argue as the proof of Corollary \ref{ssg}. For any closed point $x\in X$, let $F:=f^{-1}(x)$ and $\mathcal{J}$ be the ideal sheaf of $F$ in $X'$. Then the theorem of formal functions tells us
$$ (R^if_*\O_{X'}(E))^{\wedge}_x=\lim_{\longleftarrow}H^i(X',\O_{X'}(E)\otimes \O_{X'}/\mathcal{J}^n),$$
and the Serre duality gives us that $H^i(X',\O_{X'}(E)\otimes \O_{X'}/\mathcal{J}^n)$ is the dual of
$$ \Ext^{\dim X-i}_{X'}(\O_{X'}/\mathcal{J}^n,\omega_{X'}\otimes \O_{X'}(-E))=\Ext^{\dim X-i}_{X'}(\O_{X'}(E)\otimes \O_{X'}/\mathcal{J}^n,\omega_{X'}).$$
By the definition of local cohomology, the direct limit
$$ \lim_{\longrightarrow}\Ext^{\dim X-i}_{X'}(\O_{X'}(E)\otimes \O_{X'}/\mathcal{J}^n,\omega_{X'})$$
coincides with $H^{\dim X-i}_F(X',\omega_{X'}\otimes \O_{X'}(-E))$.

\smallskip

Using the Leray spectral sequence gives us
$$ E^{st}_2=H^s_x(X,R^tf_*(\omega_{X'}\otimes \O_{X'}(-E)))\implies H^{s+t}_F(X',\omega_{X'}\otimes \O_{X'}(-E)),$$
The relative Kawamata-Viehweg vanishing theorem gives us
$$ H^{\dim X-i}_F(X',\omega_{X'}\otimes \O_{X'}(-E))=H^{\dim X-i}_x(X,f_*(\omega_{X'}\otimes \O_{X'}(-E))).$$
In particular, since $R^if_*\O_{X'}=0$ for $1\le i<q$, $H^{\dim X-i}_x(X,f_*\omega_{X'})=0$ for $1\le i<q$.

\smallskip

By the injection $\omega_{X'}\otimes \O_{X'}(-E)\to \omega_{X'}$, we obtain an injection $f_*(\omega_{X'}\otimes \O_{X'}(-E))\to f_*\omega_{X'}$ and an exact sequence
$$ 0\to f_*(\omega_{X'}\otimes \O_{X'}(-E))\to f_*\omega_{X'}\to \mathcal{Q}\to 0$$
for some coherent sheaf $\mathcal{Q}$ on $X$. Note that the dimension of the support of $\mathcal{Q}$ is $\le \dim X-q-1$. Hence, long exact sequence of local cohomology implies that
$$ H^{\dim X-i}_x(X,f_*(\omega_{X'}\otimes \O_{X'}(-E)))=H^{\dim X-i}_x(X,f_*\omega_{X'})=0$$
for $1\le i<q$, where we used $H^{\dim X-i-1}_x(X,\mathcal{Q})=0$, which follows from dimension count. Hence $(R^if_*\O_{X'}(E))^{\wedge}_x=0$ for any closed point $x\in X$ and $R^if_*\O_{X'}(E)=0$ for $1\le i<q$.
\end{proof}

\begin{lemma}\label{b}
Let $X$ be a Cohen-Macaulay variety with $(R_q)$, and $f:X'\to X$ be any proper birational morphism from a smooth variety $X'$. Let us assume that there exists an $f$-anti-ample effective $f$-exceptional divisor $E$ on $X'$ such that $R^if_*\O_{X'}(tE)=0$ for every $1\le i<q$ and every $t\gg 0$. Then $f$ is $q$-birational.
\end{lemma}

\begin{proof}
The proof is almost verbatim to \cite[The proof of Lemma 3.3]{Kovac}.

\smallskip

By \cite[Lemma 3.8]{SchTa}, there is an injection $f_*\omega_{X'}\to \omega_X$. Moreover, by the injection $\omega_{X'}\otimes \O_{X'}(-tE)\to \omega_{X'}$, we have the following injection $f_*(\omega_{X'}\otimes \O_{X'}(-E))\to \omega_X$.

\smallskip

From this injection, there is an exact sequence
$$ 0\to f_*(\omega_{X'}\otimes \O_{X'}(-tE))\to \omega_X\to \mathcal{Q}\to 0$$
for some coherent sheaf $\mathcal{Q}$ on $X$.

\smallskip

We may take $R\mathcal{H}om_{\O_X}(-,\omega_X)$ to the exact sequence. Then the induced exact triangle is
$$ R\mathcal{H}om_{\O_{X}}(\mathcal{Q},\omega_X)\to R\mathcal{H}om_{\O_X}(\omega_X,\omega_X)\to R\mathcal{H}om_{\O_X}(f_*(\omega_{X'}\otimes \O_{X'}(-tE)),\omega_X)\to .$$

\smallskip

Now, we have
$$ \mathcal{H}^i(R\mathcal{H}om_{\O_X}(\mathcal{Q},\omega_X))=\mathcal{E}xt^i_{\O_X}(\mathcal{Q},\omega_X).$$
Moreover, by the definition of the dualizing sheaf,
$$ R\mathcal{H}om_{\O_X}(\omega_X,\omega_X)=\O_X.$$
For the third term, the relative Kawamata-Viehweg vanishing theorem gives us
$$ f_*(\omega_{X'}\otimes \O_{X'}(-tE))=Rf_*(\omega_{X'}\otimes \O_{X'}(-tE))$$
and Grothendieck duality tells us that
$$
\begin{aligned}
R\mathcal{H}om_{\O_X}(f_*(\omega_{X'}\otimes \O_{X'}(-tE)),\omega_X) &=R\mathcal{H}om_{\O_X}(Rf_*(\omega_{X'}\otimes \O_{X'}(-tE)),\omega_X)
\\ &=Rf_*R\mathcal{H}om_{\O_{X'}}(\omega_{X'}\otimes \O_{X'}(-tE),\omega_{X'})
\\ &=Rf_*\mathcal{H}om_{\O_{X'}}(\O_{X'}(-tE),\O_{X'})
\\ &=Rf_*\O_{X'}(tE).
\end{aligned}
$$

\smallskip

Hence, we have
$$ R^if_*\O_{X'}(tE)=\mathcal{E}xt^{i+1}_{\O_X}(\mathcal{Q},\omega_X).$$
Note that the codimension of the support of $\mathcal{Q}$ is the same as the codimension of the center of $f$. Thus, we can use \cite[Theorem 3.5.7]{CM}.
\end{proof}

\begin{proof}[Proof of Theorem \ref{-1application}]
This is a combination of Lemma \ref{b} and Lemma \ref{zerothapplication}.
\end{proof}

\begin{remark}
We might expect that Theorem \ref{-1application} is true for any variety $X$ with $(R_q)$ and $(S_{q+1})$. However, we cannot prove or disprove it.
\end{remark}

\begin{corollary}\label{pll}
Let $X,X'$ be any normal projective varieties, $f:X'\to X$ be any proper birational morphism, $X'$ be smooth, and $X$ be $(R_q),(S_{q+1})$. If $f$ is $q$-birational, then for any Cartier divisor $D$ on $X$ and $f$-anti-nef effective $f$-exceptional divisor $F$,
\begin{equation} \label{F}
H^i(X',\O_{X'}(f^*D+F))=H^i(X,\O_X(D))=H^i(X',\O_{X'}(f^*D))
\end{equation}
for $1\le i<q$.
\end{corollary}

\begin{proof}
The second equality is due to Proposition \ref{good}; hence, it suffices to prove the first equality only.

\smallskip

If $f$ is $q$-birational, by the following Leray spectral sequence
$$ E^{ij}_2=H^i(X,R^jf_*\O_{X'}(f^*D+F))\implies H^{i+j}(X',\O_{X'}(f^*D+F))$$
and Lemma \ref{zerothapplication}, we have the assertion.
\end{proof}

\bibliographystyle{habbvr}
\bibliography{z} 

\end{document}